\documentclass[12pt,leqno]{amsart}
\newtheorem{thm}{Theorem}[section]
\newtheorem{defi}{Definition}[section]

\newtheorem{cor}{Corollary}[section]

\newcommand{\be}{\begin{equation}}
\newcommand{\ee}{\end{equation}}
\numberwithin{equation}{section}
\newcommand{\bea}{\begin{eqnarray}}
\newcommand{\eea}{\end{eqnarray}}
\newcommand{\beb}{\begin{eqnarray*}}
\newcommand{\eeb}{\end{eqnarray*}}
\usepackage{amssymb,amsfonts,amsthm,setspace,indentfirst}
\usepackage[dvips]{graphics}
\usepackage{epsfig}
\begin{document}
\title[On Warped Product Super Generalized Recurrent Manifolds]{\bf{On Warped Product Super Generalized Recurrent Manifolds}}
\author[Absos Ali Shaikh, Haradhan Kundu and Md. Showkat Ali]{$^1$Absos Ali Shaikh, $^1$Haradhan Kundu and $^2$Md. Showkat Ali}
\date{\today}
\address{\noindent\noindent $^1$Department of Mathematics,\newline University of 
Burdwan, Golapbag,\newline Burdwan-713104,\newline West Bengal, India}
\email{aask2003@yahoo.co.in}
\email{kundu.haradhan@gmail.com}
\address{\noindent\noindent $^2$Department of Mathematics,\newline  Dhaka University,\newline  Dhaka-1000, Bangladesh.}
\email{msa317@yahoo.com}
%
%
\begin{abstract}
The object of the present paper is to obtain the characterization of a warped product semi-Riemannian manifold with a special type of recurrent like structure, called super generalized recurrent. As consequence of this result we also find out the necessary and sufficient conditions for a warped product manifold to satisfy some other recurrent like structures such as weakly generalized recurrent manifold, hyper generalized recurrent manifold etc. Finally as a support of the main result, we present an example of warped product super generalized recurrent manifold.
\end{abstract}
%
\maketitle
%
\section{\bf Introduction}\label{intro}
To generalize the notion of a manifold of constant curvature Cartan \cite{Ca26} first introduced the notion of local symmetry which can be presented as the curvature restriction $\nabla R = 0$ (i.e., the Riemann-Christoffel curvature tensor $R$ is covariantly constant). But there are many manifolds which does not bear local symmetry and hence to investigate the type of symmetry of such manifolds it is necessary to generalize the notion of local symmetry. During the last eight decades various researchers are working on this area to generalize or extend the notion of local symmetry by weakening its curvature restriction in different directions.\\
\indent Cartan \cite{Ca46} himself first gave a proper generalization of local symmetry and introduced the notion of semisymmetry, which was later classified by Szab$\acute{\mbox{o}}$ \cite{Sz82}. Then in 1983 Adam\'{o}w and Deszcz \cite{AD83} generalized the notion of semisymmetry and introduced the notion of pseudosymmetry, also known as Deszcz pseudosymmetry (see \cite{SDHJK}). On the other hand as a direct generalization of local symmetry, Chaki \cite{Ch87} introduced the notion of pseudosymmetry. We note that the interrelation between two types of pseudosymmetry is studied by Shaikh et. al. \cite{SDHJK}. In 1989 Tam\'assy and Binh \cite{TB89} generalized the Chaki's notion of pseudosymmetry and introduced the notion of weakly symmetric manifold. We note that Shaikh and his co-authors (\cite{HMS10}, \cite{SH08}-\cite{SJE08}) studied this notion of weak symmetry with various generalized curvature tensors. Again generalizing the results of Binh \cite{Binh93}, recently, Shaikh and Kundu \cite{SK12} obtained the characterization of warped product weakly symmetric manifold.\\ 
\indent Again recurrent manifold (\cite{Ruse46}, \cite{Ruse49}, \cite{Ruse49a}) is another type of generalization of local symmetry. In 1979 Dubey defined generalized recurrent manifold (briefly, $GK_n$). Recently \cite{OO12} Olszak and Olszak showed that every $GK_n$ is concircularly recurrent and every concircularly recurrent manifold is again a $K_n$ and hence every $GK_n$ is a $K_n$. Again as a generalization of $K_n$, recently, Shaikh and his coauthors introduced three new type of generalized recurrent structures together with their proper existence, namely, quasi generalized recurrent manifold (briefly, $QGK_n$) \cite{SR10}, hyper-generalized recurrent manifold (briefly, $HGK_n$) \cite{SP10} and weakly generalized recurrent manifold (briefly, $WGK_n$) \cite{SR11}. For the existence of such structures we refer the reader to see \cite{SAR13}, \cite{SRK15a}. Very recently Shaikh et. al. (\cite{SK14}, \cite{SRK15}) introduced another generalization of recurrent manifold, called, super generalized recurrent manifold  (briefly, $SGK_n$) which also generalizes the notion of $HGK_n$ as well as $WGK_n$. These kinds of generalization of recurrent structures may be called as recurrent like structures.\\
\indent The main object of the present paper is to obtain the necessary and sufficient condition for a warped product manifold to be $HGK_n$ and $WGK_n$. For this purpose we first determine the necessary and sufficient condition for a warped product to be $SGK_n$ and as consequence of this result we find out the corresponding results for $HGK_n$, $WGK_n$ and $K_n$. We know that decomposable or product manifold is a special case of warped product manifold when the warping function is identically 1. Thus we can present the characterization of a decomposable manifold with various recurrent like structures.\\
\indent The paper is organized as follows: Section 2 deals with rudimentary facts of various recurrent like structures. Section 3 is concerned with basic curvature relations of a warped product manifold. In section 4 we present our main result and, finally, in section 5 a proper example of a warped product $SGK_n$ is presented.
\section{\bf Preliminaries}\label{preli}
\indent Let $M$ be a non-flat $n$-dimensional $(n\geq 3)$ smooth manifold with semi-Riemannian metric $g$, Levi-Civita connection be $\nabla$, Riemann-Christoffel curvature tensor $R$, Ricci tensor $S$ and scalar curvature $\kappa$. In this section we define various necessary terms and curvature restricted geometric structures and for this purpose at first we consider some notations:\\
$C^{\infty}(M)$ = the algebra of all smooth functions on $M$, \\
$\chi(M)$ = the Lie algebra of all smooth vector fields on $M$, \\
$\chi^*(M)$ = the Lie algebra of all smooth 1-forms on $M$ and \\
$\mathcal T^r_k(M)$ = the space of all smooth tensor fields of type $(r,k)$ on $M$.\\
For $A,E\in \mathcal T^0_2(M)$ we have their Kulkarni-Nomizu product \cite{SK14} $A\wedge E\in \mathcal T^0_4(M)$ as
\be\label{eq2.1}
A \wedge E_{ijkl} = A_{il}E_{jk} + A_{jk}E_{il}-A_{ik}E_{jl} - A_{jl}E_{ik}.
\ee
In particular we can get $g\wedge g$, $g\wedge S$ and $S\wedge S$ as follows:
$$(g\wedge g)_{ijkl} = 2(g_{il}g_{jk}-g_{ik}g_{jl}),$$
$$(g\wedge S)_{ijkl} = g_{il}S_{jk}+S_{il}g_{jk}-g_{ik}S_{jl}-S_{ik}g_{jl} \ \mbox{ and}$$
$$(S\wedge S)_{ijkl} = 2(S_{il}S_{jk}-S_{ik}S_{jl}).$$
\begin{defi}
The manifold $M$ is said to be recurrent \cite{Wa50} if
\be\label{kn}
\nabla R = \Pi \otimes R \ \ \mbox{($\otimes$ denotes the tensor product)}
\ee
$$(\mbox{or locally, } R_{ijkl,m}=\Pi_m R_{ijkl}, \ \ \mbox{where `,' denotes the covariant derivative})$$
holds on $\{x\in M : \nabla R \neq 0 \,\,\mbox{at}\,\, x\}$ for an 1-form $\Pi\in\chi^*(M)$, called the associated 1-form of the recurrent structure. Such an $n$-dimensional manifold is denoted by $K_{n}$.
\end{defi}
\indent Again $M$ is said to be concircularly recurrent recurrent if its concircular curvature tensor $W = R - \frac{\kappa}{2n(n-1)}g\wedge g$ satisfies the condition
\be\label{ckn}
\nabla W = \Pi \otimes W
\ee
on $\{x\in M : \nabla W \neq 0 \,\,\mbox{at}\,\, x\}$ for an 1-form $\Pi\in\chi^*(M)$, called the associated 1-form.\\
The manifold $(M, g)$ is said to be generalized recurrent \cite{Du79} if it satisfies
\be\label{gkn}
\nabla R = \Pi \otimes R + \Theta\otimes g \wedge g
\ee
$$(\mbox{or locally, } R_{ijkl,m} = \Pi_m R_{ijkl} + 2 \Theta_m (g_{il}g_{jk}-g_{ik}g_{jl})$$
on $\{x\in M : \nabla R \neq \xi \otimes R \,\,\mbox{at}\,\, x\ \forall\ \xi \in \chi^*(M)\}$ for some 1-forms $\Pi$ and $\Theta$. The 1-forms $\Pi$ and $\Theta$ are called the associated 1-forms of this structure.  Such an $n$-dimensional manifold is denoted by $GK_{n}$. In \cite{OO12} Olszak and Olszak showed that every $GK_n$ satisfying \eqref{gkn} is concircularly recurrent with \eqref{ckn} and every concircularly recurrent manifold is again a $K_n$ with same associated 1-form and thus $\Theta=0$. Hence the structure $GK_n$ reduces to $K_n$. Consequently the notion of $GK_n$ does not exist.\\
\begin{defi}
The manifold $M$ is said to be quasi generalized recurrent \cite{SR10}, hyper generalized recurrent \cite{SP10}, weakly generalized recurrent \cite{SR11} and super generalized recurrent respectively if the condition
\be\label{qgkn}
\nabla R = \Pi \otimes R + \Psi \otimes \left[g\wedge (g + \eta\otimes\eta)\right],
\ee
\be\label{hgkn}
\nabla R = \Pi \otimes R + \Psi \otimes g\wedge S,
\ee
\be\label{wgkn}
\nabla R = \Pi \otimes R + \Phi \otimes S\wedge S \ \ \mbox{and}
\ee
\be\label{sgkn}
\nabla R = \Pi \otimes R + \Phi \otimes S\wedge S + \Psi \otimes g\wedge S + \Theta \otimes g\wedge g
\ee
holds respectively on $\{x\in M:\nabla R \neq \xi \otimes R \,\,\mbox{at}\,\, x\ \forall\ \xi \in \chi^*(M)\}  \subset M$ for some $\Pi$, $\Phi$, $\Psi$, $\Theta$ and $\eta \in \chi^*(M)$, called the associated 1-forms.
\end{defi}
\noindent An $n$-dimensional manifold satisfying \eqref{qgkn} is denoted by $QGK_{n}$ with $(\Pi, \Theta, \eta)$ or simply  $QGK_{n}$.\\
An $n$-dimensional manifold satisfying \eqref{hgkn} is denoted by $HGK_{n}$ with $(\Pi, \Psi)$ or simply  $HGK_{n}$.\\
An $n$-dimensional manifold satisfying \eqref{wgkn} is denoted by $WGK_{n}$ with $(\Pi, \Phi)$ or simply  $WGK_{n}$.\\
An $n$-dimensional manifold satisfying \eqref{sgkn} is denoted by $SGK_{n}$ with $(\Pi, \Phi, \Psi, \Theta)$ or simply  $SGK_{n}$.\\
\noindent In terms of local coordinates \eqref{qgkn}-\eqref{sgkn} can be respectively written as:
\be
R_{ijkl,m} = \Pi_m R_{ijkl} + \Theta_m [2(g_{il}g_{jk}-g_{ik}g_{jl}) 
					+ g_{il}\eta_j \eta_k + g_{jk}\eta_i \eta_l - g_{ik}\eta_j \eta_l-g_{jl}\eta_i \eta_k],
\ee
\be
R_{ijkl,m} = \Pi_m R_{ijkl} + \Psi_m [g_{il}S_{jk}+S_{il}g_{jk}-g_{ik}S_{jl}-S_{ik}g_{jl}],
\ee
\be
R_{ijkl,m} = \Pi_m R_{ijkl} + 2\Psi_m [S_{il}S_{jk}-S_{ik}S_{jl}] \ \ \mbox{and}
\ee
\bea\label{sgknl}
R_{ijkl,m} &=& \Pi_m R_{ijkl} + 2\Phi_m [S_{il}S_{jk}-S_{ik}S_{jl}]\\\nonumber
					 &+& \Psi_m [g_{il}S_{jk}+S_{il}g_{jk}-g_{ik}S_{jl}-S_{ik}g_{jl}]+ 2 \Theta [g_{il}g_{jk}-g_{ik}g_{jl}].
\eea
\indent It is obvious that the above structures $QGK_n$, $WGK_n$ and $HGK_n$ are all generalization of $K_n$ but all three exists independently (see \cite{SAR13}, \cite{SRK15a}). We also mention that $SGK_n$ is a proper generalization of $WGK_n$ and $HGK_n$ (see Section \ref{exam}).
\begin{defi} 
Let $D\in \mathcal T^0_4(M)$ and $A, E, F \in \mathcal T^0_2(M)$. Then $M$ is said to be Roter type manifold (briefly, $RT_n$) with $(D;A,E)$ (\cite{Desz03}, \cite{Desz03a}) and generalized Roter type manifold (briefly, $GRT_n$) with $(D;A,E,F)$ (\cite{SDHJK}, \cite{SKgrt}, \cite{SKgrtw}) respectively if
$$D = N_1 A\wedge A - N_2 A\wedge E - N_3 E\wedge E \ \mbox{ and}$$
$$D = L_1 A\wedge A - L_2 A\wedge E - L_3 E\wedge E - L_4 A\wedge F - L_5 E\wedge F - L_6 F\wedge F$$
respectively, for some $N_i, L_j \in C^{\infty}(M)$, $1\le i\le 3$ and $1\le j\le 6$.
\end{defi}
%
\section{\bf Warped Product Manifold}\label{warp}
In 1957 Kru$\breve{\mbox{c}}$kovi$\breve{\mbox{c}}$ \cite{Kr57} initiated the study of semi-decomposable manifolds which were latter named as warped product manifolds by Bishop and O'Neill \cite{BO69}.
Let $(\overline M, \overline g)$ and $(\widetilde M, \widetilde g)$ be two semi-Riemannian manifolds of dimension $p$ and $(n-p)$ respectively ($1\leq p < n$), and $f$ is a positive smooth  function on $\overline M$. Then the warped product $M= \overline M\times_{f}\widetilde M$ is the product manifold $\overline M\times \widetilde M$  of dimension $n$ endowed with the metric
\be\label{warpf}
g=\pi^*(\overline g) + (f\circ\pi) \sigma^* (\widetilde g),
\ee
where $\pi:M\rightarrow\overline M$ and $\sigma:M\rightarrow\widetilde M$ are the natural projections. Then the local components of the metric $g$ are given by
\begin{eqnarray}\label{eq3.1}
g_{ij}=\left\{\begin{array}{lll}
&\overline g_{ij}&\ \ \ \ \mbox{for} \ i = a \ \mbox{and} \ j = b,\\
&f \widetilde g_{ij}&\ \ \ \ \mbox{for $i = \alpha$ and $j = \beta$,}\\
&0&\ \ \ \ \mbox{otherwise.}\\
\end{array}\right.
\end{eqnarray}
Here $a,b \in \left\{1,2,...,p\right\}$ and $\alpha, \beta \in \left\{p+1,p+2,...,n\right\}$. We note that throughout the paper we consider $a,b,c,...\in \{1,2, ..., p\}$ and $\alpha,\beta,\gamma,...\in \{p+1,p+2,...,n\}$ and $i,j,k,...\in \{1,2,...,n\}$. Here $\overline M$ is called the base, $\widetilde M$ is called the fiber and $f$ is called warping function of the warped product $M = \overline M \times_f \widetilde M$. It may be mentioned that the warped product metric $g$ can be taken as (see \cite{Geba93}, \cite{Kr57}, \cite{Zh14})
$$g=\pi^*(\overline g) + (f\circ\pi)^2 \sigma^* (\widetilde g).$$
However throughout the paper we will consider the warped product metric given in \eqref{warpf}.
Again we assume that, when $\Omega$ is a quantity formed with respect to $g$, we denote by $\overline \Omega$ and $\widetilde \Omega$, the similar quantities formed with respect to $\overline g$ and $\widetilde g$ respectively. By straightforward calculation one can easily calculate the local components of $\Omega$ in terms of $\overline \Omega$ and $\widetilde \Omega$ for $\Omega = \nabla, R, S$ and $\kappa$ and obtain the following:\\
The non-zero local components of Levi-Civita connection $\nabla$ of $M$ are given by
\be\label{eq2.2}
\Gamma^a_{bc}=\overline{\Gamma}^a_{bc},\,\,\,\, \Gamma^\alpha_{\beta \gamma}=\widetilde{\Gamma}^\alpha_{\beta \gamma},\,\,\,\,\,\,\,\Gamma^a_{\beta \gamma}=-\frac{1}{2}\overline{g}^{ab}f_{b} \widetilde{g}_{\beta \gamma},\,\,\ \ \Gamma^\alpha_{ a \beta }=\frac{1}{2f}f_{a}\delta^{\alpha}_{\beta},
\ee
where $f_{a}=\partial_{a} f=\frac{\partial f}{\partial x^{a}}$.\\ 
The local components of the Riemann-Christoffel curvature tensor $R$ and Ricci tensor $S$ of $M$ which may not vanish identically are the following:
\be\label{eq2.3}
R_{abcd} = \overline{R}_{abcd},\,\,\,\, R_{a\alpha b\beta}=f T_{ab}\widetilde{g}_{\alpha \beta},\,\,\,R_{\alpha \beta \gamma \delta} = f\widetilde{R}_{\alpha \beta \gamma \delta} - f^2 P \widetilde{G}_{\alpha \beta \gamma \delta}, 
\ee
\be\label{eq2.4}
S_{ab}=\overline{S}_{ab}-(n-p)T_{ab},\,\,\,\, S_{\alpha \beta}=\widetilde{ S}_{\alpha \beta} + Q \widetilde{g}_{\alpha \beta},
\ee
where $G_{ijkl} = \frac{1}{2}g\wedge g_{ijkl} = g_{il}g_{jk}-g_{ik}g_{jl}$ are the components of Gaussian curvature and
$$T_{ab} = -\frac{1}{2f}(\nabla_b f_a - \frac{1}{2f}f_a f_b), \ \ \ \ tr(T) = g^{ab}T_{ab},$$
$$Q = f((n-p-1)P -tr(T)), \ \ \ \ P = \frac{1}{4f^2}g^{ab}f_a f_b.$$
The scalar curvature $\kappa$ of $M$ is given by
\be\label{eq2.5}
\kappa=\bar{\kappa}+\frac{\tilde{\kappa}}{f}-(n-p)[(n-p-1)P - 2 \; tr(T)].  
\ee
Again the non-zero local components of $\nabla R$ are given by \cite{Ho04}:
\bea\label{eq2.6}
&&\left\{\begin{array}{l}
(i) R_{abcd,e} = \overline R_{abcd,e},\\
(ii) R_{a\alpha b\beta,e} = f T_{ab,e} \widetilde g_{\alpha \beta},\\
(iii) R_{\alpha\beta\gamma\delta,e} = -f_e \widetilde R_{\alpha\beta\gamma\delta} + f^2 P_e \widetilde G_{\alpha\beta\gamma\delta},\\
(iv) R_{\alpha\beta\gamma\delta,\epsilon} = f \widetilde R_{\alpha\beta\gamma\delta,\epsilon},\\
(v) R_{\alpha\beta\gamma d,\epsilon} = -\frac{f_d}{2} \widetilde R_{\alpha\beta\gamma\epsilon} + \frac{f^2}{2} P_d \widetilde G_{\alpha\beta\gamma\epsilon},\\
(vi) R_{abc\delta,\epsilon} = \frac{1}{2}\widetilde g_{\epsilon\delta}(f_a T_{bc}-f_ b T_{ac}) + \frac{1}{2}f^d R_{abcd} \widetilde g_{\epsilon\delta}.\\
\end{array}\right.
\eea
The non-zero components of $(g\wedge g)$, $(g\wedge S)$ and $(S\wedge S)$ are given by
\bea\label{gg}
\left\{
\begin{array}{l}
(i) (g\wedge g)_{abcd} = (\overline g \wedge \overline g)_{abcd},\\
(ii) (g\wedge g)_{a\alpha b \beta} = -2f \overline g_{ab}\widetilde g_{\alpha\beta},\\
(iii) (g\wedge g)_{\alpha\beta\gamma\delta} = f^2 (\widetilde g \wedge \widetilde g)_{\alpha\beta\gamma\delta},
\end{array}
\right.
\eea
\bea\label{gs}
\left\{
\begin{array}{l}
(i) (g\wedge S)_{abcd} = (\overline g\wedge \overline S)_{abcd} - (n-p) (\overline g \wedge T)_{abcd},\\
(ii) (g\wedge S)_{a\alpha b \beta} = -\overline g_{ab}(\widetilde S_{\alpha\beta}+Q\widetilde g_{\alpha\beta})-f \widetilde g_{\alpha\beta}[\overline S_{ab}-(n-p)T_{ab}],\\
(iii) (g\wedge S)_{\alpha\beta\gamma\delta} = f (\widetilde g\wedge \widetilde S)_{\alpha\beta\gamma\delta} + f Q \widetilde (\widetilde g \wedge \widetilde g)_{\alpha\beta\gamma\delta},
\end{array}
\right.
\eea
\bea\label{ss}
\left\{
\begin{array}{l}
(i) (S\wedge S)_{abcd} = (\overline S\wedge \overline S)_{abcd} - (n-p) (\overline S \wedge T)_{abcd} + (n-p)^2 (T \wedge T)_{abcd},\\
(ii) (S\wedge S)_{a\alpha b \beta} = -2(\widetilde S_{\alpha\beta}+Q\widetilde g_{\alpha\beta})[\overline S_{ab}-(n-p)T_{ab}],\\
(iii) (S\wedge S)_{\alpha\beta\gamma\delta} = (\widetilde S\wedge \widetilde S)_{\alpha\beta\gamma\delta} + Q \widetilde (\widetilde g\wedge \widetilde S)_{\alpha\beta\gamma\delta} + Q^2 (\widetilde g \wedge \widetilde g)_{\alpha\beta\gamma\delta}.
\end{array}
\right.
\eea
For detailed information about the local components of various tensors on a warped product manifold we refer the reader to see \cite{On83}, \cite{SK12}, \cite{SKgrtw} and also references therein.
%
\section{\bf Warped Product $SGK_n$}\label{main}
\begin{thm}\label{thm4.1}
Let $M^n = \overline M^p \times_f \widetilde M^{n-p}$ be a warped product manifold. Then $M$ is a $SGK_n$ with $(\Pi,\Phi,\Psi,\Theta)$ if and only if the following conditions hold simultaneously
\bea\label{eq4.1}
\left\{
\begin{array}{l}
(i) \overline\nabla\overline R = \overline \Pi \otimes \overline R + \overline \Phi \otimes \overline S \wedge \overline S 
																		+ \overline \Psi \otimes \overline g\wedge \overline S 
																		+ \overline \Theta \otimes \overline g \wedge \overline g\\
								\hspace{0.5in} - 2(n-p)\overline \Phi \otimes \overline S \wedge T + (n-p)^2\overline \Phi \otimes T\wedge T
																		- (n-p)\overline \Psi \otimes \overline g \wedge T,\\
(ii) -\widetilde \Pi \otimes \overline R = \widetilde \Phi \otimes \overline S \wedge \overline S
						+\widetilde \Psi \otimes \overline g \wedge \overline S + \widetilde \Theta \otimes \overline g \wedge \overline g\\
 \hspace{0.9in}	- 2(n-p)\widetilde \Phi \otimes \overline S \wedge T + (n-p)^2\widetilde \Phi \otimes T \wedge T
 					- (n-p)\widetilde \Psi \otimes \overline g \wedge T,
\end{array}
\right.
\eea
\bea\label{eq4.2}
\left\{
\begin{array}{l}
(i) -(df + f \overline \Pi)\otimes \widetilde R = \overline \Phi \otimes \widetilde S \wedge \widetilde S 
																							+ (2 Q \overline \Phi + f\overline \Psi) \otimes \widetilde g \wedge \widetilde S\\
															\hspace{1.35in}+ (-\frac{1}{2} f^2(P \overline \Pi + dP) + Q^2 \overline \Phi + f Q \overline \Psi
																							+ f^2\overline \Theta) \otimes \widetilde g \wedge \widetilde g,\\
(ii) f \widetilde\nabla\widetilde R = f \widetilde \Pi \otimes \widetilde R 
																						+ \widetilde \Phi \otimes \widetilde S \wedge \widetilde S
																						+ (Q \widetilde \Phi + f\widetilde \Psi) \otimes \widetilde g \wedge \widetilde S\\ 
															\hspace{0.7in}+ (-\frac{1}{2} f^2 P \widetilde \Pi + Q^2 \widetilde \Phi + f Q \widetilde \Psi 
																								+ f^2 \widetilde \Theta) \otimes \widetilde g \wedge \widetilde g,\\
\end{array}
\right.
\eea
\bea\label{eq4.3}
\left\{
\begin{array}{l}
(i) \left[2 \overline \Phi\otimes(\overline S - (n-p)T) + \overline \Psi\otimes \overline g\right]\otimes \widetilde S =\\
					\hspace{0.1in} -[f(\overline \nabla T - \overline \Pi\otimes T)+ (2Q \overline \Phi+ f\overline \Psi)
					\otimes (\overline S - (n-p)T) + (Q\overline \Psi + 2f \overline \Theta )\otimes \overline g] \otimes \widetilde g,\\
(ii) \left[2 \widetilde \Phi\otimes(\overline S -(n-p)T) +\widetilde \Psi\otimes \overline g\right]\otimes\widetilde S=\\
							\hspace{0.2in} -[-f\widetilde \Pi\otimes T+ (2Q \widetilde \Phi+ f\widetilde \Psi)\otimes 
								(\overline S - (n-p)T) + (Q\widetilde \Psi + 2f \widetilde \Theta )\otimes \overline g] \otimes \widetilde g,\\
\end{array}
\right.
\eea
\bea\label{eq4.4}
\left\{
\begin{array}{l}
(i) f^d R_{abcd} = -(f_a T_{bc}-f_b T_{ac}) \ \ \mbox{and}\\
(ii) df\otimes \widetilde R = f^2 \Theta P\otimes \widetilde G.\\
\end{array}
\right.
\eea
\end{thm}
\noindent {\bf Proof:} First suppose that $M$ is $SGK_n$. Then in terms of local coordinates the defining condition can be written as
\be\label{sgk}
R_{ijkl,m} = \Pi_m R_{ijkl} + \Phi_m (S\wedge S)_{ijkl} + \Psi_m (g\wedge S)_{ijkl} + \Theta_m (g\wedge g)_{ijkl}.
\ee
Putting 
\beb
\left\{
\begin{array}{l}
(i) i=a,j=b,k=c,l=d,m=e \ \ \mbox{and}\\
(ii) i=a,j=b,k=c,l=d,m=\epsilon\\
\end{array}
\right.
\eeb
respectively in \eqref{sgk} and then in view of \eqref{eq2.3}-\eqref{ss} it is easy to check that \eqref{eq4.1} holds. Similarly putting 
\beb
\left\{
\begin{array}{l}
(iii) i=\alpha,j=\beta,k=\gamma,l=\delta,m=e;\\
(iv) i=\alpha,j=\beta,k=\gamma,l=\delta,m=\epsilon;\\
(v) i=a,j=\alpha,k=b,l=\beta,m=e;\\
(vi) i=a,j=\alpha,k=b,l=\beta,m=\epsilon;\\
(vii) i=a,j=b,k=c,l=\alpha,m=\epsilon \ \ \mbox{and}\\
(viii) i=\alpha,j=\beta,k=\gamma,l=a,m=\epsilon\\
\end{array}
\right.
\eeb
respectively in \eqref{sgk} and then in view of \eqref{eq2.3}-\eqref{ss} we get \eqref{eq4.2}-\eqref{eq4.4} respectively. The converse part is obvious. This proves the theorem.\\
\indent From Theorem \ref{thm4.1}, it follows that the nature of base and fiber of a warped product $SGK_n$ is given by the following:
\begin{cor}\label{cor4.1}
Let $M^n = \overline M^p \times_f \widetilde M^{n-p}$ be a warped product $SGK_n$with $(\Pi,\Phi,\Psi,\Theta)$. Then\\
(i) $\overline M$ is a $SGK_p$ if $T$ can be expressed as a linear combination of $\overline S$ and $\overline g$.\\
(ii) $\overline M$ is a $GRT_p$ with $(\overline R; \overline g, \overline S, T)$ on 
$\{x \in M: \widetilde \Pi_x \ne 0\}$.\\
(iii) $\widetilde M$ is a $SGK_{n-p}$ with $(\widetilde \Pi,\frac{1}{f}\widetilde \Phi,\frac{Q}{f}\widetilde \Phi+\widetilde \Psi,-\frac{1}{2} f P \widetilde \Pi + \frac{Q^2}{f} \widetilde \Phi + Q \widetilde \Psi + f \widetilde \Theta)$.\\
(iv) $\widetilde M$ is a $RT_{n-p}$ with $(\widetilde R; \widetilde g, \widetilde S)$ on 
$\{x \in M: (df + f \overline \Pi)_x \ne 0\}$.\\
(v) $\widetilde M$ satisfies Einstein metric condition on 
$\{x \in M: (2(\overline \kappa - (n-p)tr(T)) \overline \Phi + p\overline \Psi)_x \ne 0\}$ $\cup$
$\{x \in M: (2(\overline \kappa - (n-p)tr(T)) \widetilde \Phi + p\widetilde \Psi)_x \ne 0\} =$
$\{x \in M: (2(\overline \kappa - (n-p)tr(T)) \Phi + p \Psi)_x \ne 0\}$.\\
(vi) $\widetilde M$ is of constant curvature on $\{x \in M: df_x \ne 0\}$.
\end{cor}
From Theorem \ref{thm4.1}, the characterization of a decomposable or product $SGK_n$ is given by the following:
\begin{cor}\label{cor4.2}
Let $M^n = \overline M^p \times \widetilde M^{n-p}$ be a product manifold. Then $M$ is a $SGK_n$ with $(\Pi,\Phi,\Psi,\Theta)$ if and only if\\
$1. (i)$ $\overline\nabla\overline R = \overline \Pi \otimes \overline R + \overline \Phi \otimes \overline S \wedge \overline S 
																		+ \overline \Psi \otimes \overline g\wedge \overline S 
																		+ \overline \Theta \otimes \overline g \wedge \overline g$,\\
\indent (ii) $-\widetilde \Pi \otimes \overline R = \widetilde \Phi \otimes \overline S \wedge \overline S
						+\widetilde \Psi \otimes \overline g \wedge \overline S + \widetilde \Theta \otimes \overline g \wedge \overline g$,\\
$2. (i)$ $\widetilde\nabla\widetilde R = \widetilde \Pi \otimes \widetilde R+\widetilde\Phi\otimes\widetilde S\wedge \widetilde S 
																		+ \widetilde \Psi \otimes \widetilde g\wedge \widetilde S 
																		+ \widetilde \Theta \otimes \widetilde g \wedge \widetilde g$,\\
\indent (ii) $-\overline \Pi \otimes \widetilde R = \overline \Phi \otimes \widetilde S \wedge \widetilde S
						+\overline \Psi \otimes \widetilde g \wedge \widetilde S +\overline \Theta \otimes \widetilde g \wedge \widetilde g$,\\
$3. (i)$ $\left[2\overline \Phi\otimes\overline S + \overline \Psi\otimes \overline g\right]\otimes \widetilde S =
					-[\overline \Psi\otimes \overline S + 2\overline \Theta\otimes \overline g] \otimes \widetilde g,$\\
\indent (ii) $\left[2\widetilde \Phi\otimes\widetilde S + \widetilde \Psi\otimes \widetilde g\right]\otimes \overline S =
					-[\widetilde \Psi\otimes \widetilde S + 2\widetilde \Theta\otimes \widetilde g] \otimes \overline g.$
\end{cor}
From the above, the nature of each factor of a product $SGK_n$ is given by the following:
\begin{cor}\label{cor4.3}
Let $M^n = \overline M^p \times \widetilde M^{n-p}$ be a product $SGK_n$ with $(\Pi,\Phi,\Psi,\Theta)$. Then\\
(i) base and fiber are both super generalized recurrent manifolds.\\
(ii) $\overline M$ is a $RT_p$ with $(\overline R; \overline g, \overline S)$ on 
$\{x \in M: \widetilde \Pi_x \ne 0\}$.\\
(iii) $\widetilde M$ is a $RT_{n-p}$ with $(\widetilde R; \widetilde g, \widetilde S)$ on 
$\{x \in M: \overline \Pi_x \ne 0\}$.\\
(iv) $\overline M$ satisfies Einstein metric condition on 
$\{x \in M: (2\widetilde \kappa \widetilde \Phi + (n-p)\widetilde \Psi)_x \ne 0\}$.\\
(v) $\widetilde M$ satisfies Einstein metric condition on 
$\{x \in M: (2\overline \kappa \overline \Phi + p\overline \Psi)_x \ne 0\}$.
\end{cor}
\begin{cor}\label{cor4.4}\cite{SK12}
Let $M^n = \overline M^p \times_f \widetilde M^{n-p}$ be a warped product manifold. Then $M$ is a recurrent manifold with
$$\nabla R = \Pi \otimes R$$
if and only if the following conditions hold simultaneously\\
$1. (i)$ $\overline \nabla \overline R = \overline \Pi \otimes \overline R,$
\indent (ii) $\widetilde \Pi \otimes \overline R = 0,$\\
$2. (i)$ $-(df + f \overline \Pi)\otimes \widetilde R =\frac{1}{2}f^2(P \overline \Pi - dP)\otimes\widetilde g\wedge \widetilde g,$
\indent (ii) $\widetilde\nabla\widetilde R = \widetilde \Pi \otimes \widetilde R$ and $P \widetilde \Pi  =0,$\\
$3. (i)$ $\overline \nabla T = \overline \Pi\otimes T,$
\indent (ii) $\widetilde \Pi\otimes T = 0,$\\
$4. (i)$ $f^d R_{abcd} = -(f_a T_{bc}-f_b T_{ac})$ and (ii) $df\otimes \widetilde R = f^2 dP\otimes \widetilde G$.
\end{cor}
\begin{cor}\label{cor4.5}
Let $M^n = \overline M^p \times_f \widetilde M^{n-p}$ be a warped product recurrent manifold satisfying
$\nabla R = \Pi \otimes R$. Then\\
$(i)$ $\overline M$ and $\widetilde M$ are both recurrent. Also $T$ is recurrent with associated 1-form $\Pi$.\\
$(ii)$ $\overline R = 0$, $T =0$ and $P =0$ on the set $\{x\in M:\widetilde \Pi \ne 0\}$.\\
$(iii)$ $\widetilde M$ is of constant curvature on $\{x\in M:df_x \ne 0\} \cup \{x\in M:(df + f \overline \Pi)_x \ne 0\}$.
\end{cor}
\begin{cor}\label{cor4.6}
Let $M^n = \overline M^p \times \widetilde M^{n-p}$ be a product manifold. Then $M$ is recurrent satisfying
$$\nabla R = \Pi \otimes R$$
if and only if\\
$(i)$ $\overline \nabla \overline R = \overline \Pi \otimes \overline R,$
\indent $\widetilde \Pi \otimes \overline R = 0,$
\indent $(ii)$ $\overline \Pi\otimes \widetilde R = 0,$
\indent $\widetilde\nabla\widetilde R = \widetilde \Pi \otimes \widetilde R.$
\end{cor}
\begin{cor}\label{cor4.7}
Let $M^n = \overline M^p \times_f \widetilde M^{n-p}$ be a warped product manifold. Then $M$ is a $HGK_n$ with $(\Pi,\Psi)$
if and only if the following conditions hold simultaneously\\
$1. (i)$ $\overline \nabla \overline R = \overline \Pi \otimes \overline R 
																				+ \overline \Psi \otimes \overline g \wedge \overline S 
																				- (n-p)\overline \Psi \otimes \overline g \wedge T,$\\
\indent (ii) $-\widetilde \Pi \otimes \overline R = \widetilde \Psi \otimes \overline g \wedge \overline S
																									- (n-p)\widetilde \Psi \otimes \overline g \wedge T,$\\
$2. (i)$ $-(df + f \overline \Pi)\otimes \widetilde R = f\overline \Psi \otimes \widetilde g \wedge \widetilde S 
																	+ (-\frac{1}{2} f^2(P \overline \Pi + dP) + f Q \overline \Psi) 
																	\otimes \widetilde g \wedge \widetilde g,$\\
\indent (ii) $f \widetilde\nabla\widetilde R = f \widetilde \Pi \otimes \widetilde R 
																					+ f\widetilde \Psi \otimes \widetilde g \wedge \widetilde S 
																		+ (-\frac{1}{2} f^2 P \widetilde \Pi + f Q \widetilde \Psi)
																		\otimes \widetilde g \wedge \widetilde g,$\\
$3. (i)$ $\overline \Psi\otimes \overline g\otimes \widetilde S =
											-[f(\overline \nabla T - \overline \Pi\otimes T) + f\overline \Psi\otimes (\overline S - (n-p)T) + Q\overline \Psi\otimes \overline g] \otimes \widetilde g,$\\
\indent (ii) $\widetilde \Psi\otimes \overline g \otimes \widetilde S = -[-f\widetilde \Pi\otimes T + f\widetilde \Psi \otimes (\overline S - (n-p)T) + Q\widetilde \Psi \otimes \overline g] \otimes \widetilde g,$\\
$4. (i)$ $f^d R_{abcd} = -(f_a T_{bc}-f_b T_{ac})$ and\\
\indent (ii) $df\otimes \widetilde R = f^2 dP\otimes \widetilde G$.
\end{cor}
\begin{cor}\label{cor4.8}
Let $M^n = \overline M^p \times_f \widetilde M^{n-p}$ be a warped product $HGK_n$ with $(\Pi,\Psi)$. Then\\
$(i)$ $\overline M$ is hyper generalized recurrent manifold if $T$ and $S$ are linearly dependent.\\
$(ii)$ $\widetilde M$ is a $SGK_{n-p}$ with $(\widetilde \Pi,0,\widetilde \Psi,-\frac{1}{2} f P \widetilde \Pi + Q \widetilde \Psi)$.\\
$(iii)$ $\widetilde M$ is of vanishing conformal curvature tensor on $\{x\in M:(df + f \overline \Pi)_x \ne 0\}$.\\
$(iv)$ $\widetilde M$ satisfies Einstein metric condition on $\{x\in M:\Psi_x \ne 0\}$ and $\widetilde M$ is of constant curvature on $\{x\in M:df_x \ne 0\}$.
\end{cor}
\noindent \textbf{Proof:} Results of (i), (ii) and (iv) are obvious from Corollary \ref{cor4.7}. By virtue of 2.(i) of Corollary \ref{cor4.7}, on $\{x\in M:(df + f \overline \Pi)_x \ne 0\}$, $\widetilde R$ can be exressed as a linear combination of $\widetilde g \wedge \widetilde S$ and $\widetilde g \wedge \widetilde g$. Hence in view of Corollary 6.1 of \cite{SK14}, the conformal curvature tensor of $\widetilde M$ vanishes on this set, which proves (iii).
\begin{cor}\label{cor4.9}
Let $M^n = \overline M^p \times \widetilde M^{n-p}$ be a product manifold. Then $M$ is a $HGK_n$ with $(\Pi,\Psi)$ if and only if the following conditions hold simultaneously\\
$1. (i)$ $\overline \nabla \overline R = \overline \Pi \otimes \overline R 
																				+ \overline \Psi \otimes \overline g \wedge \overline S,$ \ \ \ 
		(ii) $\widetilde \Pi \otimes \overline R + \widetilde \Psi \otimes \overline g \wedge \overline S = 0,$\\
$2. (i)$ $\overline \Pi\otimes \widetilde R + \overline \Psi \otimes \widetilde g \wedge \widetilde S = 0,$ \ \ \ \ \ 
		(ii) $\widetilde\nabla\widetilde R = \widetilde \Pi \otimes \widetilde R 
																		+ \widetilde \Psi \otimes \widetilde g \wedge \widetilde S,$\\
$3. (i)$ $\overline \Psi\otimes \overline g\otimes \widetilde S = -\overline \Psi\otimes \overline S \otimes \widetilde g,$ \ \ \ 
		(ii) $\widetilde \Psi\otimes \overline g \otimes \widetilde S = -\widetilde \Psi \otimes \overline S \otimes \widetilde g$.
\end{cor}
\begin{cor}\label{cor4.10}
Let $M^n = \overline M^p \times_f \widetilde M^{n-p}$ be a warped product manifold. Then $M$ is a $WGK_n$ with $(\Pi,\Phi)$
if and only if the following conditions hold simultaneously\\
$1. (i)$ $\overline \nabla \overline R = \overline \Pi \otimes \overline R 
																			+ \overline \Phi \otimes \overline S \wedge \overline S - 2(n-p)\overline \Phi 
																			\otimes (\overline S\wedge T) + (n-p)^2\overline \Phi \otimes T\wedge T,$\\
\indent (ii) $-\widetilde \Pi \otimes \overline R = \widetilde \Phi \otimes \overline S \wedge \overline S  
																				- 2(n-p)\widetilde \Phi\otimes (\overline S\wedge T) 
																				+ (n-p)^2\widetilde \Phi \otimes T\wedge T,$\\
$2. (i)$ $-(df + f \overline \Pi)\otimes \widetilde R = \overline \Phi \otimes \widetilde S \wedge \widetilde S 
																											+ 2 Q \overline \Phi \otimes \widetilde g \wedge \widetilde S 
																											+ (-\frac{1}{2}f^2(P \overline \Pi - dP) 
																											+ Q^2 \overline \Phi) \otimes \widetilde g \wedge \widetilde g,$\\
\indent (ii) $f \widetilde\nabla\widetilde R = f \widetilde \Pi \otimes \widetilde R 
																					+ \widetilde \Phi \otimes \widetilde S \wedge \widetilde S 
																		+ Q \widetilde \Phi\otimes \widetilde g \wedge \widetilde S 
													+ (-\frac{1}{2} f^2 P \widetilde \Pi + Q^2 \widetilde \Phi) \otimes \widetilde g \wedge \widetilde g,$\\
$3. (i)$ $2 \overline \Phi\otimes(\overline S - (n-p)T) \otimes \widetilde S =
			-[f(\overline \nabla T - \overline \Pi\otimes T)+ 2Q\overline \Phi\otimes (\overline S - (n-p)T)] \otimes \widetilde g,$\\
\indent (ii) $2 \widetilde \Phi\otimes(\overline S - (n-p)T)\otimes \widetilde S =
							-[-f\widetilde \Pi\otimes T+ 2 Q\widetilde \Phi\otimes (\widetilde S - (n-p)T)] \otimes \widetilde g,$\\
$4. (i)$ $f^d R_{abcd} = -(f_a T_{bc}-f_b T_{ac})$ and\\
\indent (ii) $df\otimes \widetilde R = f^2 dP\otimes \widetilde G$.
\end{cor}
\begin{cor}\label{cor4.11}
Let $M^n = \overline M^p \times_f \widetilde M^{n-p}$ be a warped product $WGK_n$ with $(\Pi,\Phi)$. Then\\
$(i)$ $\overline M$ is a $WGK_{p}$ if $T$ and $S$ are linearly dependent.\\
$(ii)$ $\widetilde M$ is a $SGK_{n-p}$ with $(\widetilde \Pi,\Phi,\frac{Q}{f}\widetilde \Psi,-\frac{1}{2} f P \widetilde \Pi + \frac{Q^2}{f} \widetilde \Phi)$.\\
$(iii)$ $\widetilde M$ is a $RT_{n-p}$ with $(\widetilde R;\widetilde g,\widetilde S)$ on $\{x\in M:(df + f \overline \Pi)_x \ne 0\}$.\\
$(iv)$ $\widetilde M$ satisfies Einstein metric condition on
$\{x \in M: ((\overline \kappa - (n-p)tr(T)) \Phi)_x \ne 0\}$
and is of constant curvature on $\{x\in M:df_x \ne 0\}$.
\end{cor}
\begin{cor}\label{cor4.12}
Let $M^n = \overline M^p \times \widetilde M^{n-p}$ be a product manifold. Then $M$ is a $WGK_n$ with $(\Pi,\Phi)$ if and only if the following conditions hold simultaneously\\
$1. (i)$ $\overline \nabla \overline R = \overline \Pi \otimes \overline R 
																			+ \overline \Phi \otimes \overline S \wedge \overline S,$ \ \ \ 
		(ii) $\widetilde \Pi \otimes \overline R + \widetilde \Phi \otimes \overline S \wedge \overline S = 0,$\\
$2. (i)$ $\overline \Pi \otimes \widetilde R + \overline \Phi \otimes \widetilde S \wedge \widetilde S = 0,$ \ \ \ 
		(ii) $\widetilde\nabla\widetilde R = \widetilde \Pi \otimes \widetilde R 
																			+ \widetilde \Phi \otimes \widetilde S \wedge \widetilde S,$\\ 
$3. (i)$ $\overline \Phi\otimes \overline S \otimes \widetilde S = 0$, \ \ \ 
		(ii) $\widetilde \Phi\otimes \overline S \otimes \widetilde S = 0$.
\end{cor}
\section{\bf An example of warped product $SGK_4$}\label{exam}
\textbf{Example 1:} Consider the warped product $M = \overline M\times_f \widetilde M$, where $\overline M$ is a 3-dimensional manifold equipped with the metric
$$\overline{ds}^2 = e^{x^2}(dx^1)^2+e^{x^1}(dx^2)^2+(dx^3)^2$$
in local coordinates $(x^1, x^2, x^3)$ and $\widetilde M$ is an open interval of $\mathbb R$ with local coordinate $x^4$ and the warping function $f=e^{x^3}$. The non-zero components of the Riemann-Christoffel curvature tensor $\overline R$ and Ricci tensor $\overline S$ of $\overline M$ are given by
$$ \overline R_{1212}=-\frac{1}{4} \left(e^{x^1}+e^{x^2}\right), \ \ 
\overline S_{11}=\frac{1}{4} \left(1+e^{x^2-x^2}\right), \ \overline S_{22}=\frac{1}{4} \left(e^{x^1-x^2}\right).$$
Again the non-zero components of $\overline \nabla \overline R$ are
$$\overline R_{1212,1}=\frac{e^{x^2}}{4},~~ \overline R_{1212,2}=\frac{e^{x^1}}{4}.$$
If we consider the 1-form $\overline \Pi = \left(-\frac{e^{x^2}}{e^{x^1}+e^{x^2}},-\frac{e^{x^1}}{e^{x^1}+e^{x^2}},0\right)$, then we can easily check that the manifold $\overline M$ is recurrent satisfying  $\overline\nabla \overline R = \overline \Pi\otimes \overline R$.\\
Now the metric of $M = \overline M\times_f \widetilde M$ is given by
$$ds^2 = e^{x^2}(dx^1)^2+e^{x^1}(dx^2)^2+(dx^3)^2+e^{x^3}(dx^4)^2.$$
The non-zero local components of the Riemann-Christoffel curvature tensor $R$ and Ricci tensor $S$ of $M$ are given by
$$ R_{1212}=-\frac{1}{4} \left(e^{x^1}+e^{x^2}\right), ~~ 
R_{3434}=-\frac{e^{x^3}}{4},$$
$$S_{11}=\frac{1}{4} \left(e^{x^2-x^1}+1\right), ~~ 
S_{22}=\frac{1}{4} \left(e^{x^1-x^2}+1\right), ~~ 
S_{33}=\frac{1}{4}, ~~ 
S_{44}=\frac{e^{x^3}}{4}.$$
Again the non-zero local components of $\nabla R$ are given by
$$ R_{1212,1}=\frac{e^{x^2}}{4}, ~~ 
R_{1212,2}=\frac{e^{x^1}}{4}.$$
Again the non-zero components of $g\wedge g$, $g\wedge S$ and $S\wedge S$ given as follows:
$$ g\wedge g_{1212}=-2e^{x^1+x^2}, ~~ 
g\wedge g_{1313}=-2e^{x^2}, ~~ 
g\wedge g_{1414}=-2e^{x^2+x^3}, ~~ 
g\wedge g_{2323}=-2e^{x^1},$$ 
$$g\wedge g_{2424}=-2e^{x^1+x^3}, ~~ 
g\wedge g_{3434}=-2e^{x^3},$$
$$ g\wedge S_{1212}=\frac{1}{2} \left(-e^{x^1}-e^{x^2}\right), ~~ 
g\wedge S_{1313}=\frac{1}{4} \left(-e^{x^2-x^1} \left(e^{x^1}+1\right)-1\right),$$
$$ g\wedge S_{1414}=-\frac{1}{4} e^{x^3-x^1} \left(e^{x^2} \left(e^{x^1}+1\right)+e^{x^1}\right), ~~ 
g\wedge S_{2323}=\frac{1}{4} \left(-e^{x^1-x^2} \left(e^{x^2}+1\right)-1\right),$$
$$ g\wedge S_{2424}=-\frac{1}{4} e^{x^3-x^2} \left(e^{x^2} \left(e^{x^1}+1\right)+e^{x^1}\right), ~~ 
g\wedge S_{3434}=-\frac{e^{x^3}}{2},$$
$$ S\wedge S_{1212}=\frac{1}{4} \left(-\cosh \left(x^1-x^2\right)-1\right), ~~ 
S\wedge S_{1313}=\frac{1}{8} \left(-e^{x^2-x^1}-1\right),$$
$$ S\wedge S_{1414}=-\frac{1}{8} e^{x^3} \left(e^{x^2-x^1}+1\right), ~~ 
S\wedge S_{2323}=\frac{1}{8} \left(-e^{x^1-x^2}-1\right),$$
$$ S\wedge S_{2424}=-\frac{1}{8} e^{x^3} \left(e^{x^1-x^2}+1\right), ~~ 
S\wedge S_{3434}=-\frac{e^{x^3}}{8}.$$
If we consider the 1-forms $\Pi$, $\Phi$, $\Psi$ and $\Theta$ as:
\be
\Pi_{i}=\left\{\begin{array}{ccc}
&\frac{ \Psi_1 e^{x^1} \left (e^{x^2} - 2 \right) + 2  \Psi_1 \cosh (x^1 - x^2) - (2  \Psi_1 + 1) e^{x^2} + 2  \Psi_1} {2 \left (e^{x^1} + e^{x^2} \right)}&\mbox{for} \ i = 1\\

&\frac { \Psi_2 \left (e^{x^1} - 2 \right) e^{x^2} + 2  \Psi_2 \cosh (x^1 - x^2) - (2  \Psi_2 + 1) e^{x^1} + 2  \Psi_2} {2 \left (e^{x^1} + e^{x^2} \right)}&\mbox{for} \ i = 2\\

&\frac { \Psi_3 e^{-x^1 - x^2} \left (-e^{x^1 + x^2} + e^{x^1} + e^{x^2} \right)^2} {2 \left (e^{x^1} + e^{x^2} \right)}&\mbox{for} \ i = 3\\

&\frac { \Psi_4 e^{-x^1 - x^2} \left (-e^{x^1 + x^2} + e^{x^1} + e^{x^2} \right)^2} {2 \left (e^{x^1} + e^{x^2} \right)}&\mbox{for} \ i = 4,\\
\end{array}\right.
\ee
\be
\Phi_{i}=\left\{\begin{array}{ccc}
&\frac {\left ( \Psi_1 \left (-e^{x^1 + x^2} \right) + 2  \Psi_1 \cosh (x^1 -x^2) + 2  \Psi_1 + e^{x^2} \right)} {-2 \cosh (x^1 - x^2) + \sinh (x^1) + \cosh (x^1) + \sinh (x^2) + \cosh (x^2) - 2}&\mbox{for} \ i = 1\\

&e^{x^1} \left (\frac { \Psi_2 e^{x^1} + 1} {e^{x^1} + e^{x^2}} -  \Psi_2 + \frac {1} {\left (e^{x^1} - 1 \right) e^{x^2} - e^{x^1}} \right) - \Psi_2&\mbox{for} \ i = 2\\

&-\frac {\Psi_3 \left (e^{x^1 + x^2} + e^{x^1} + e^{x^2} \right)} {e^{x^1} + e^{x^2}}&\mbox{for} \ i = 3\\

&-\frac {\Psi_4 \left (e^{x^1 + x^2} + e^{x^1} + e^{x^2} \right)} {e^{x^1} + e^{x^2}}&\mbox{for} \ i = 4,\\
\end{array}\right.
\ee
\be
\Theta_{i}=\left\{\begin{array}{ccc}
&\frac {-\Psi_1 e^{x^1 - x^2} + 2  \Psi_1 e^{x^2} \sinh (x^1) - 2  \Psi_1 + e^{x^2}} {16 \left (-e^{x^1 + x^2} + e^{x^1} + e^{x^2} \right)}&\mbox{for} \ i = 1\\

&\frac {1} {16} \left(- \Psi_2 e^{-x^1} -  \Psi_2 e^{-x^2} -  \Psi_2 + \frac {e^{x^1}} {-e^{x^1 + x^2} + e^{x^1} + e^{x^2}} \right)&\mbox{for} \ i = 2\\

&-\frac {1} {16}  \Psi_3 e^{-x^1 - x^2} \left (e^{x^1 + x^2} + e^{x^1} + e^{x^2} \right)&\mbox{for} \ i = 3\\

&-\frac {1} {16}  \Psi_4 e^{-x^1 - x^2} \left (e^{x^1 + x^2} + e^{x^1} + e^{x^2} \right)&\mbox{for} \ i = 4\\
\end{array}\right.
\ee
then we can check that $M$ is a $SGK_4$ with $(\Pi, \Phi, \Psi, \Theta)$, which is neither a $HGK_4$ nor a $WGK_4$.\\
\noindent
\textbf{Acknowledgment:} 
The second named author gratefully acknowledges to CSIR, New Delhi (File No. 09/025 (0194)/2010-EMR-I) for the financial assistance. All the algebraic computations of Section 5 are performed by a program in Wolfram Mathematica.


\end{document}